# LINGUISTIC PARADOXES AND TAUTOLOGIES

by Florentin Smarandache
University of New Mexico
Gallup, NM 87301, USA

**Abstract.**
Classes of linguistic paradoxes are introduced with examples and explanations. They are part of the author's work on the Paradoxist Philosophy based on mathematical logic.
The general cases exposed below are modeled on the English language structure in a rigid way. In order to find nice particular examples of such paradoxes one grammatically adjusts the sentences.

**A. Introduction**:
Let <N>, <V>, <A> be some noun, verb, and attribute respectively, and <Non-N>, <Non-V>, <Non-A> respectively their antonyms. For example, if <A> is <small> then <Non-A> is <big> or <large>, etc.
Also let <N'>, <N''>, etc. represent synonyms of <N> or even <N>, and so <V'>, <V''>, etc.
or <A'>, <A''>, etc.
Let <NV> represent a noun-ed verb, and <NV'> a synonym, etc.

**B. Classes of Paradoxes**:
Then, one defines the following classes of linguistic paradoxes and semi-paradoxes:

1. **All is <A>, the <Non-A> too.**
Examples:
All is possible, the impossible too.
All is real, the unreal too.
All is justice, the injustice too.
All life is complex, the simple life too.
All people are actors, the non-actors too.
All can be happy, even the unhappy.
It's so near, but yet so far away.
All is weird, the natural too.
All is joyful, the sorrow too.

2. **<Non-N> is a better <N>.**
   **<Non-A> is a better <A>.**
   **<Non-V> is a better <V>.**
Examples:
Not to speak is sometimes a better speech.
Not to complain is a better complain.
Unattractive is sometimes better than attractive.
Slow is sometimes better than fast.

No government is a better government.
A non-ruler is a better ruler.
No news is good news.
Not to stare is sometimes better to look.
Not to love is a better love.
Not to move is sometimes a better move.
Impoliteness is a better politeness.
Not to hear is better than not listening.
No reaction is sometimes the best reaction.
Not to show kindness is a better kindness [welfare].
She is better than herself.
No fight is a better fight [i.e. to fight by non-violent means].

**3. Only <N> is truly a <Non-N>.**
**Only <A> is truly a <Non-A>.**
Examples:
Only a rumor is truly a gossip.
Only a fiction is truly a fact.
Only normal is truly not normal.
Nobody is truly a 'somebody'.
Only fiction is truly real.
The friend is the most dangerous hidden enemy.
Only you are truly not you [=you act strangely].
Only mercy can be truly merciless.
If you spit at the sky, it will fall in your face.
Only gentleness is truly wild.

**4. This is so <A>, that it looks <Non-A>.**
Examples:
This is so true, that it looks false!
This is so ripe, that it looks spoilt.
This is so friendly, that it looks hostile.
He seemed so trustworthy, that he looks untrustworthy.
This is so fake, it looks real!
This is so proper, that it looks improper.
This is so beautiful, that it looks unreal.
This is so simple, that it looks difficult.
The story was so real, that it looked fiction.
Can't see the trees for the forest.

**5. There is some <N> which is <A> and <Non-A> at the same time.**
Examples:
There are events which are good and bad at the same time.
There are laws which are good and bad at the same time.
There are some news which are real and wrong at the same time.
There are some insects which are helpful and dangerous at the same time. [like the spider]
There are men who are handsome and ugly at the same time.
There are classes that are fun and boring at the same time.
There are some ministers which are believers and mis-believers at the same time.
There are moments that are sweet and sour.
There are games which are challenged and not competitive at the same time.
Food which are simultaneously hot and cold.
The game was exciting, yet boring [because we were losing].
People are smart and foolish at the same time [i.e., smart at something, and foolish at other thing].

**6. There is some <N> which <V> and really <Non-V> at the same time.**
   Examples:
There are people who trick and do not really trick at the same time.
There are some people who play and don't play at the same time.
Some of life's experiences are punishments and rewards at the same time.
Exercise is exhausting but also invigorating.
There are children who listen and do not really listen at the same time.
There are teachers who teach and don't teach at the same time.
There are people who spell and misspell at the same time.
Nice and rough people concomitantly.
Politicians who lie and tell the truth all the time!

**7. To <V>, even when <Non-V>.**
   Examples:
A sage thinks even when he doesn't think.
I exist when I don't exist.
A clown is funny even when he isn't being funny.
To die of thirst surrounded by water.  [saltwater]
To be a poet and not know it.
A mother worries even when she doesn't worry.
To believe even when you don't believe.
Is matching even when not-matching.
I sleep even when I am awake.
Always running around.
To dream, even when not sleeping.

**8. This <N> is enough <Non-N>.**
   Examples:
This silence is enough noise.
This vacation from work is hard work.
The superiority brings enough inferiority.  [=listlessness]
This day is my night.
This diary is enough non-diary.
This sleep is enough awake.
This sweet truthfulness is enough sarcasm.
This table of four is enough for six people.
I had enough.
This job is enough recreation [when enjoying the job].

**9.   <Non-V> sometimes means <V>.**
  Examples:
Not to speak sometimes means to speak.
Not to touch sometimes means to touch.
The preserve peace sometimes means going to war.
To destroy life (as in viruses) sometimes means to preserve life.
Not to listen sometimes means to listen.
Two feet forward sometimes means standing still.
Not to litter sometimes means to litter.
Speeding is sometimes not speeding [in case of emergency].
Not to show anger is sometimes to show anger.

**10.   <N> without <N>.**
  Examples:
Hell without hell.
The style without style.
The rule applied:  there were no rule!
Our culture is our lack of culture.
Live without living.
Some people are so afraid of death, that they do not live.
Work without work.
Can't live with them, can't live without them.
Death without death [for a Christian dying is going on to eternal life].
Guilt without guilt [sometimes is guilty but doesn't feel guilty].

**11.   a) <N> inside/within the <Non-N>.**
  Examples:
Movement inside the immobility.
Silence within the noise.
Slavery within the freedom.
Loneliness within a crowd.
A circle within a circle.
The wrestling ring inside a squared section.
To find wealth in poverty [i.e., happiness and love].

       **b) <Non-N> in the <N>.**
  Examples:
Immobility inside the movement.
Noise inside the silence.
The eye of the storm.
Government.  Bureaucracies.
Inequality inside the equality.
Single inside the marriage.
Anger inside the happiness.
Warmth in the cold.
Cold in the heat.
Laughing without being happy.
Has not gotten anywhere.

Poverty in wealth [no poverty or love in a wealth family].

**12. The <A> of the <Non-A>.**
  Examples:
The shadow of the light.
Music of silence.
Relaxing of exercise effect.
The restrictions of the free.
Life through death.
The sound/loudness of the silence.
I can see the light at the of the tunnel.
The slave of freedom [someone who couldn't give up his freedom, not even in marriage].

**13. <V> what one <Non-V>.**
  Examples:
To see what one can't see.
To hear what one can't hear.
To taste what one can't taste.
To accept what on can't understand.
To say what one can't say [to tell a secret].
To wait patiently when one doesn't know how to wait.
To breath what one can't breath.
To feel what one can't feel.
To appreciate what one dis-appreciates.
To believe what one can't believe [faith].
To smell what one can't smell.

**14. Let's <V> by <Non-V>.**
  Examples:
Let's strike by not striking [=Japanese strike].
Let's talk by not talking; [means to think].
To vote by not voting at all.
To help someone by not helping [using experience as a teacher].
Let's justify by not justifying.
Let's win by not winning.
Let's strip by not stripping [to make bare or clear].
Let's fight by not fighting [Ghandi's Motto].

**15. <N> of the <Non-N>.**
  Examples:
The benefits we get from non-benefits.
The smoke we got from non-smokers.
The rewards we get from hard work.
The service we get from non-service.
The god that comes from bad.
The pleasure we get from the pain.

**16. <Non-A> is <A>.**

Examples:
The bad is good [because makes you try harder].
The good is bad [because doesn't leave any room for improvement].
Work is a blessing.
The poor is spiritually rich.
Sometimes ugly is beauty [because beauty is in the eyes of the beholden].
You have to kiss a lot of frogs before you find a prince.
Hurt is healing.
"There is no absolute" is an absolute.
Not to commit any error is an error.

**17. A <Non-N> <N>.**
Examples:
A positive negative [which means: a failure enforces you to do better]
A sad happiness.
An impossible possibility.
Genuine imitation leather.
A loud whisper.
A beautiful disaster [which means beauty can be found anywhere].
A bitter sweet.
A harsh gentleness [a gentleness that is very firm with you].
A guiltless sinner [someone who doesn't regret sinning].

**18. Everything has an <A> and a <Non-A>.**
Examples:
Everything has a sense and a non-sense.
Everything has a truthful side and a wrong side.
Everything has a beginning and an ending.
Everything has a birth and a death.
Everything has its time and a non-time.
Everything has an appearance and a non-appearance.
Everything around you resolves and also dissolves.
Everybody has a good side and a bad side.
Everyone has a right and a wrong.

**19. <V> what <Non-V>.**
Examples:
To be what you are not.
One needs what one doesn't need.
Expect the unexpected!
Culture exists by its non-existence.
No matter how rich we are, we never make enough money.
One purchases what one doesn't purchase.
To work when we are not working.
To die might mean to live for ever [= for an artist].
One wants because one doesn't want [sometimes one wants something only because someone else likes it].

**Linguistic Tautologies:**

   A tautology is a redundacy, a pleonasm, a needless repetition of an idea, according to the "Webster's New World Dictionary", Third College Edition, 1988.
However, the following classes of tautologies - using repetition - go to a deeper meaning, and even changes the sense.  A double assertation reverses to a negation.
One also may play with the synonyms.

   **20.   Mirror semi-paradox:**
         **<N> of the <N'>.**
         **<N> of the <N'> of the <N''>...**
   Examples:
Best of the best.
Worst of the worst of the worst.
Mother of the mother of the mother... [the maternal grandmother]
Follower of the followers.
The rows of the rows [lots of rows.].

   **21.   This is not an <N>, this is an <N'>.**
   Examples:
This is not a teacher, this is a professor.
This is not a car, this is a Wolswagen.
This is not a truck, this is a Chevy.
This is not noise, this is music.
This is not music, this is noise.
This is not a cedar tree, this is a <gad>  [gad = Navajo name for cedar tree].
This is not me, this is I.
This is not a sword, this is a saber.
This is not a problem, this is an exercise [= easier].
Practice makes you practice.
This is not a girl, this is Katie.
This is not a horse, this is a pony.

   **22.   <N> is not enough <N'>.**
         **<A> is not enough <A'>.**
   Examples:
Sufficient is not enough sufficient [which means: to do more than "sufficient"].
Punishment is not enough punishment.
Health is not enough wealth.
Clean is not enough clean.
Studying is not enough studying [which means to do more than just getting by, i.e. to do research].
Extravagant is not enough extravagant.
Time is not enough time.
The more you have, the more you want.
Attention is not enough attention [some people need action too].

**23. More <A> than <A'>.**
  Examples:
Better than better [=perfection].
Worst than worst [=evil].
Sweeter than sweeter [=honey].
More life than life [=spirituality].
More depressed than depressed.
Foster than foster.
More beautiful than pretty.
More ugly than ugly [really ugly].
Smarter than smart [like a genius].

**24. How <A> is an <A'> <N>?**
  Examples:
How democratic is a so called democratic society?
How republican is a so called republic society?
How civilized is a so called civilized person?
How free is a free country?
How commanding is a so called commanding officer?
How Pop Culture is a so called Pop Culture?
How strong is a strong man?
How lone is a lone ranger? [not very, he has tanto].

**25. No <A> is really <A'>.**
  Examples:
No friend is really a friend [s/he betrays you when you don't even expect!].
No luck is really a luck.
No original is really original.
No husband is really a husband [you learn to depend on yourself!].
No tomboy is really a tomboy [girl considered boyish].
No work is really less work.
No true Marxist is really a true Marxist [they contradict their own beliefs].
No magic is really magic [all is only a trick].

**26. I would rather prefer <A>, than <A'>.**
  Examples:
I would rather prefer pretty, than prettier.
I want that, not that.
I would rather prefer this, than this.
I would rather be old, than *old*.
I would rather prefer great than big.
I would rather be crazy than crazy [crazy like foolish, than crazy like insane].

**27. More <A> than <A'>.**
  Examples:

Prettier than pretty.
More real than real.
More advantage than advantage.
More help than help.
More smiles than smiles [she didn't psychically smile, but there were smiles written all over her face].
More cries than cries.
More meters than kilometers.
Make everyday a rainbow day.
He earns more than himself.
More suspicious than suspected.

   28.  <V> those who <V'> you.
  Examples:
Ignore those who ignore you.
Criticize those who criticize you.
Defend those who defend you.

   29.  <V>, because <V'>.
  Examples:
I want because I want.
I think because I think.
I hear because I listen.
I see because I look.
I need because I need.
I know because I know.
I live because I live.
I believe what is unbelievable [faith].
I am happy because I am happy [there is no reason for my happiness].

   30.  <V> the <NV'>.
  Examples:
I hate the haters (therefore I hate myself!).
I envy the enviers (therefore I envy myself).
I am strange to strangers.
I cheat the cheaters (therefore I cheat myself).
I lie to liars (therefore I lie to myself).
I kick the kickers (therefore I kick myself).
I love the lovers.

  **Exercises for readers**:
Try to construct a general scheme - using <N>, <V>, <A>, etc. notations as above - and then give particular cases for each of the following paradoxes or semi-paradoxes:

- Dream the impossible dream.
- It is not a question of what we are, but more of who we are.

- Only a small dream/output is really a big dream/output.
- One vote is enough to make a difference; and yet one vote isn't often enough to make a difference.
- Good times come and good times go, but memories last for lifetime.
- The dark and light of infinity.
- The sense we get from non-sense.
- Think before you think.
- Sometimes less is more.
- Enjoy life today, tomorrow may never come.
- Make it happen, by making it happen.
- Less is more.
- My shoes are cleaner than my feet.
- No matter how hard it seems, it will get easier.
- See the things as they are not [see their hidden spot].
- Quitters never win, and winners never quit.
- My needs exist for needs.
- Bad things happen for a good reason.

Look at this Funny Law example:

**A Paradoxist Government**:
   Suppose you have two cows.  Then the government kills them and milks you!

**A poem**:

Sometimes in life we see but do not have sight
             or don't see what we should see
                we hear but do not listen
                we speak but do not communicate
                we live but do not know how to live
                we love but do not love
And then we die but we've already have been dead

   The list of such invented linguistic paradoxes can be indefinitely extended.  It is specific to each language, and it is based on language expressions and types of sentence and phrase constructions and structures.
One can also play with antonymic/synonymic adverbs, prepositions, etc. to construct other categories of linguistic paradoxes.

# 24.1.119, 114, February 1997.